\documentclass[12pt]{article}
\usepackage{amsmath,amsfonts, amssymb}
\usepackage{graphicx}
\usepackage{longtable}
\usepackage{url}
\usepackage{tocloft}
\makeatletter
 
  \@addtoreset{equation}{section}
\makeatother
\def\beqar {\begin{eqnarray}}
\def\eeqar {\end{eqnarray}}
\def\beq {\begin{equation}}
\def\eeq {\end{equation}}

\def\S{{\cal S}}

\def\al{\alpha}

\def\del{\delta}
\def\Del{\Delta}

\def\Ga{\Gamma}

\def\ep{\epsilon}
\def\la{\lambda}

\def\om{\omega}

\def\d{\partial}

\def\diag{{\rm diag}}

\def\det{{\rm det}}
\def\cp{\mathbb{CP}}
\def\C{\mathbb{C}}

\def\Z{\mathbb{Z}}

\begin{document}

\begin{titlepage}
\null\vspace{-62pt} \pagestyle{empty}
\begin{center}
\vspace{1.0truein}

{\Large\bf Elements of Aomoto's generalized hypergeometric \\
\vspace{.35cm}
functions and a novel perspective on Gauss' \\
\vspace{.35cm}
hypergeometric differential equation} \\

\vspace{1.0in} {\sc Yasuhiro Abe} \\
\vskip .12in {\it Cereja Technology Co., Ltd.\\
3-11-15 UEDA-Bldg. 4F, Iidabashi   \\
Chiyoda-ku, Tokyo 102-0072, Japan }\\
\vskip .07in {\tt abe@cereja.co.jp}\\
\vspace{1.3in}

\centerline{\large\bf Abstract}
\end{center}
We review Aomoto's generalized hypergeometric functions on Grassmannian spaces $Gr(k+1,n+1)$.
Particularly, we clarify integral representations of the
generalized hypergeometric functions in terms of twisted homology and cohomology.
With an example of the $Gr(2,4)$ case, we consider in detail Gauss' original hypergeometric
functions in Aomoto's framework. This leads us to present a new systematic
description of Gauss' hypergeometric differential equation
in a form of a first order Fuchsian differential equation.

\end{titlepage}
\pagestyle{plain} \setcounter{page}{2} 

\tableofcontents


\section{Introduction}

Studies of physical phenomena such as scattering amplitudes
(see, {\it e.g.}, \cite{ArkaniHamed:2012nw}-\cite{Ferro:2015grk})
and quantum Hall effects \cite{Balli:2014pqa} in Grassmannian spaces $Gr(k+1,n+1)$
have been attentively carried out in recent years.
This has revived an interest in a purely mathematical subject, {\it i.e.},
generalized hypergeometric functions on Grassmannian spaces, which were
introduced and developed by Gelfand \cite{Gelfand:1986} and independently
by Aomoto \cite{Aomoto:1994bk} many years ago.
One of the main goals of this chapter is to present a clear and systematic
review on these particular topics in mathematical physics.
Particularly, we clarify integral representations of
Aomoto's generalized hypergeometric functions in terms of twisted homology and cohomology.
As the simplest example, we consider in detail Gauss' original hypergeometric
functions in Aomoto's framework so as to familiarize ourselves to
the concept of twisted homology and cohomology.

This chapter is organized as follows. In the next section we review
some formal results of Aomoto's generalized hypergeometric functions on $Gr(k+1,n+1)$,
based on Japanese textbooks \cite{Aomoto:1994bk,Haraoka:2002bk}.
We present a review in a pedagogical fashion since these results are
not familiar enough to many physicists and mathematicians.
In section 3 we consider a particular case $Gr (2, n+1)$
and present its general formulation.
In section 4 we further study the case of $Gr(2,4)$
which reduces to Gauss' hypergeometric function.
Imposing permutation invariance among branch points, we
obtain new realizations of the hypergeometric differential equation
in a form of a first order Fuchsian differential equation.
This is an original result first reported in \cite{Abe:2015ucn}
and this chapter is to be based mainly on some part of \cite{Abe:2015ucn}.

\section{Elements of Aomoto's generalized hypergeometric functions}
\subsection{Definition}

Let $Z$ be a $(k+1) \times (n+1)$ matrix
\beq
    Z = \left(
          \begin{array}{ccccc}
            z_{00} & z_{01} & z_{02} & \cdots & z_{0n} \\
            z_{10} & z_{11} & z_{12} & \cdots & z_{1n} \\
            \vdots & \vdots & \vdots &  & \vdots \\
            z_{k0} & z_{k1} & z_{k2} & \cdots & z_{kn} \\
          \end{array}
        \right)
    \label{1-1}
\eeq
where $k < n$ and the matrix elements are complex, $z_{ij} \in \C$
($0 \le i \le k \, ;~  0 \le j \le n $).
A function of $Z$, which we denote $F(Z)$, is defined as {\it
a generalized hypergeometric function on Grassmannian space $Gr (k+1, n+1)$
}
when it satisfies the following relations:
\beqar
    \sum_{j = 0}^{n} z_{ij} \frac{\d F}{\d z_{pj}} &=& - \del_{ip} F
    ~~~ ( 0 \le i, p \le k )
    \label{1-2a} \\
    \sum_{i = 0}^{k} z_{ij} \frac{\d F}{\d z_{ij}} &=&  \al_{j} F
    ~~~ ( 0 \le j \le n )
    \label{1-2b} \\
    \frac{\d^2 F}{\d z_{ip}\d z_{jq}} &=& \frac{\d^2 F}{\d z_{iq}\d z_{jp}}
    ~~~ ( 0 \le i, j \le k \,  ; ~ 0 \le p, q \le n  )
    \label{1-2c}
\eeqar
where the parameters $\al_j$ obey the non-integer conditions
\beqar
    \al_j & \not\in & \Z ~~~~~ ( 0 \le j \le n )
    \label{1-3a} \\
    \sum_{j = 0}^{n} \al_j & =& - (k + 1)
    \label{1-3b}
\eeqar

\subsection{Integral representation of $F(Z)$ and twisted cohomology}

The essence of Aomoto's generalized hypergeometric function
\cite{Aomoto:1994bk} is that, by use of the so-called
twisted de Rham cohomology,\footnote{
The {\it twisted} de Rham cohomology is a version of
the ordinary de Rham cohomology into which multivalued functions,
such as $\Phi$ in (\ref{1-5a}), are incorporated.
For mathematical rigor on this, see Section 2 in \cite{Aomoto:1994bk}.
}
$F(Z)$ can be written in a form of integral:
\beq
    F ( Z ) \, = \, \int_\Delta \Phi \om
    \label{1-4}
\eeq
where
\beqar
    \Phi &=& \prod_{j = 0}^{n} l_j (\tau )^{\al_j}
    \label{1-5a} \\
    l_j (\tau ) &=& \tau_0 z_{0j} + \tau_1 z_{1j}
    + \cdots + \tau_k z_{k j} ~~~~ ( 0 \le j \le n)
    \label{1-5b} \\
    \om &=& \sum_{i = 0}^{k} ( -1 )^i \tau_i
    d \tau_0 \wedge d \tau_1 \wedge \cdots \wedge
    d \tau_{i-1} \wedge d \tau_{i + 1} \wedge \cdots \wedge d \tau_k
    \label{1-5c}
\eeqar
The complex variables $\tau = ( \tau_0 , \tau_1 , \cdots , \tau_k )$
are homogeneous coordinates of the complex projective space $\cp^k$,
{\it i.e.}, $\C^{k+1 } - \{ 0, 0, \cdots , 0 \}$.
The multivalued function $\Phi$ is then defined in a space
\beq
    X \, = \, \cp^{k} - \bigcup_{j = 0}^{n} {\cal H}_j
    \label{1-6}
\eeq
where
\beq
    {\cal H}_j \, = \, \{  \tau \in \cp^k  \, ; ~ l_j (\tau ) = 0 \}
    \label{1-7}
\eeq

We now consider the meaning of the integral path $\Delta$.
Since the integrand $\Phi \om$ is a multivalued
$k$-form, simple choice of $\Delta$ as
a $k$-chain on $X$ is not enough.
{\it Upon the choice of $\Delta$, we need to implicitly
specify branches of $\Phi$ on $\Delta$ as well, otherwise
we can not properly define the integral.}
In what follows we assume these implicit conditions.

Before considering further properties of $\Delta$,
we here notice that $\om$ has an ambiguity
in the evaluation of the integral (\ref{1-4}).
Suppose $\al$ is an arbitrary $(k -1)$-form defined in $X$.
Then an integral over the exact $k$-form $d( \Phi \al )$ vanishes:
\beq
    0 \, = \, \int_\Delta d ( \Phi \al ) \, = \,
    \int_\Delta \Phi \left(
    d \al + \frac{d \Phi}{\Phi} \wedge \al
    \right) \, = \, \int_\Delta \Phi \nabla \al
    \label{1-8}
\eeq
where $\nabla$ can be interpreted as a covariant (exterior) derivative
\beq
    \nabla  \, = \, d  + d \log \Phi \wedge \, = \,
    d + \sum_{j = 0}^{n} \al_j \frac{ d l_j}{l_j} \wedge
    \label{1-9}
\eeq
This means that $\om^\prime = \om + \nabla \al$
is equivalent to $\om$ in the definition of the integral (\ref{1-4}).
Namely, $\om$ and $\om^\prime$ form an equivalent class,
$\om \sim \om^\prime$. This equivalent class is called the
cohomology class.

To study this cohomology class, we consider the differential equation
\beq
    \nabla f \, = \, d f + \sum_{j = 0}^{n} \al_j \frac{d l_j}{l_j} f \, =\, 0
    \label{1-10}
\eeq
General solutions are locally determined by
\beq
    f \, = \, \la \, \prod_{j = 0}^{n} l_{j} (\tau )^{- \al_j} ~~~~~~
    ( \la \in \C^{\times} )
    \label{1-11}
\eeq
These local solutions are thus basically given by $1 / \Phi$.
The idea of locality is essential since even if $1 / \Phi$ is multivalued
within a local patch it can be treated as a single-valued function.
Analytic continuation of these solutions forms a fundamental homotopy group
of a closed path in $X$ (or $1/X$ to be precise but it can
be regarded as $X$ by flipping the non-integer powers $\al_j$ in (\ref{1-5a})).
The representation of this fundamental group is called the {\it monodromy}
representation.
The monodromy representation determines the local system
of the differential equation (\ref{1-10}).
The general solution $f$ or $1 / \Phi$ gives a rank-1 local system in
this sense\footnote{
It is `rank-1' because each factor $l_j ( \tau )$ in the local solutions
(\ref{1-11}) is first order in the elements of $\tau$.
}.
We denote this rank-1 local system by ${\cal L}$.
The above cohomology class is then defined as an element
of the $k$-th cohomology group of $X$ over ${\cal L}$,
{\it i.e.},
\beq
    [ \om ] \in H^k ( X , {\cal L} )
    \label{1-12}
\eeq
This cohomology group $H^k ( X , {\cal L} )$ is also
called {\it twisted} cohomology group.

\subsection{Twisted homology and twisted cycles}

Having defined the cohomology group $H^k ( X , {\cal L} )$, we
can now define the dual of it, {\it i.e.}, the $k$-th
homology group $H_k ( X , {\cal L}^{\vee})$, known as the twisted
homology group, where ${\cal L}^{\vee}$ is the rank-1 dual local system
given by $\Phi$.
A differential equation corresponding to ${\cal L}^{\vee}$
can be written as
\beq
    \nabla^{\vee} g \, = \, d g -  \sum_{j = 0}^{n} \al_j \frac{d l_j}{l_j}
    g \, =\, 0
    \label{1-13}
\eeq
We can easily check that the general solutions are given by $\Phi$:
\beq
    g \, = \, \la \, \prod_{j = 0}^{n} l_{j} (\tau )^{ \al_j}
    \, = \, \la \Phi ~~~~~~
    ( \la \in \C^{\times} )
    \label{1-14}
\eeq
As before, an element of $H_k ( X , {\cal L}^{\vee})$
gives an equivalent class called a homology class.

In the following, we show that the integral path $\Del$
forms an equivalent class and see that it coincides with
the above homology class.
Applying Stokes' theorem to  (\ref{1-8}), we find
\beq
    0  \, = \, \int_\Delta \Phi \nabla \al \, = \, \int_{\d \Delta} \Phi \al
    \label{1-15}
\eeq
where $\al$ is an arbitrary $(k -1 )$-form as before.
The boundary operator $\d$ is in principle
determined from $\Phi$ (with information on branches).
Denoting $C_p ( X , {\cal L}^{\vee})$
a $p$-dimensional chain group on $X$ over ${\cal L}^{\vee}$,
we can express the boundary operator as
$\d \, : \, C_p ( X , {\cal L}^{\vee})  \longrightarrow
C_{p-1} ( X , {\cal L}^{\vee})$.
Since the relation (\ref{1-15}) holds for an arbitrary $\al$,
we find that the $k$-chain $\Del$ vanishes by the action of $\d$:
\beq
    \d \Delta \, = \, 0
    \label{1-16}
\eeq
The $k$-chain $\Del$ satisfying above is generically called the $k$-cycle.
In the current framework it is also  called the {\it twisted cycle}.
Since the boundary operator satisfies $\d^2 = 0$, the $k$-cycle
has a redundancy in it. Namely,
$\Del^\prime = \Del + \d C_{(+1)}$
also becomes the $k$-cycle where $C_{(+1)}$ is
an arbitrary $( k + 1 )$-chain or an element of $C_{k+1} ( X , {\cal L}^{\vee})$.
Thus $\Del$ and $\Del^\prime$ form an equivalent class,
$\Del \sim \Del^\prime$, and this is exactly the
homology class defined by $H_k ( X , {\cal L}^{\vee} )$, {\it i.e.},
\beq
    [ \Delta ] \in H_k ( X , {\cal L}^{\vee} )
    \label{1-17}
\eeq

To summarize, the generalized hypergeometric
function (\ref{1-4}) is determined by the following bilinear form
\beqar
    H_k ( X , {\cal L}^{\vee} ) \times H^k ( X , {\cal L} )
    & \longrightarrow & \C
    \label{1-18} \\
    \left( [ \Delta ], [ \om ] \right) & \longrightarrow &
    \int_{\Del} \Phi \om
    \label{1-19}
\eeqar

\subsection{Differential equations of $F(Z)$}

The condition $l_j (\tau ) = 0$ in (\ref{1-7}) defines
a hyperplane in $(k+1)$-dimensional spaces.
To avoid redundancy in configuration of hyperplanes,
we assume the set of hyperplanes are non-degenerate,
that is, we consider the hyperplanes in {\it general position}.
This can be realized by demanding that any $(k+1)$-dimensional
minor determinants of the $(k+1) \times (n+1)$ matrix $Z$ are nonzero.
We then redefine $X$ in (\ref{1-6}) as
\beq
    X =  \{
    Z \in Mat_{k+1 , n+1} ( \C ) | \mbox{ {\rm any $(k+1)$-dim
    minor determinants of $Z$ are nonzero}}
    \}
    \label{1-20}
\eeq
In what follows we implicitly demand this condition in $Z$.
The configuration of $n+1$ hyperplanes in $\cp^k$ is determined by
this matrix $Z$.

Apart from the concept of hyperplanes, we can also interpret that
the above $Z$ provides $n+1$ {\it distinct points} in $\cp^k$.
Since a homogeneous coordinate of $\cp^k$ is given by
$\C^{k+1} - \{0, 0, \cdots, 0 \}$, we can consider each of
the $n+1$ column vectors of $Z$ as a point in $\cp^k$;
the $j$-th column representing the $j$-th homogeneous coordinates
of $\cp^k$ ($j = 0,1, \cdots , n$).

The scale transformation, under which
the $\cp^k$ homogeneous coordinates are invariant,
is realized by an action of $H_{n+1} = \{ \diag (h_0 , h_1 , \cdots h_n )
| h_j \in \C^\times \}$ from right on $Z$.
The general linear transformation of the
homogeneous coordinates, on the other hand,
can be realized by an action of $GL(k+1, \C)$ from left.
These transformations are then given by
\beqar
    {\mbox {\rm Linear transformation:}}&~& Z \rightarrow Z^\prime = g Z
    \label{1-21a} \\
    {\mbox {\rm Scale transformation:}}&~& Z \rightarrow Z^\prime = Z h
    \label{1-21b}
\eeqar
where $g \in GL(k+1, \C)$ and $h \in H_{n+1}$.
Under these transformations the integral $F(Z)$ in (\ref{1-4}) behaves as
\beqar
    F (g Z ) &=& (\det g)^{-1} F(Z)
    \label{1-22a}\\
    F ( Z h ) &=& F (Z) \prod_{j = 0}^{n} h_{j}^{\al_j}
    \label{1-22b}
\eeqar

We now briefly show that the above relations lead to the
defining equations of the generalized hypergeometric functions
in (\ref{1-2a}) and (\ref{1-2b}), respectively.
Let ${\bf 1}_n$ be the $n$-dimensional identity matrix ${\bf 1}_n
= \diag ( 1,1, \cdots ,1)$, and $E_{ij}^{(n)}$ be an
$n \times n$ matrix in which
only the $(i,j)$-element is 1 and the others are zero.
We consider $g$ in a particular form of
\beq
    g = {\bf 1}_{k+1} + \ep E^{(k+1)}_{pi}
    \label{1-23}
\eeq
where $\ep$ is a parameter.
Then $gZ$ remains the same as $Z$ except the $p$-th row
which is replaced by $(z_{p0} + \ep z_{i0}, z_{p1} + \ep z_{i1}, \cdots,
z_{pn} + \ep z_{in})$.
Then the derivative of $F(gZ)$ with respect to $\ep$ is expressed as
\beq
    \frac{\d}{\d \ep} F (gZ) \, = \, \sum_{j=0}^{n} z_{ij}
    \frac{\d }{\d z_{pj}} F(gZ)
    \label{1-24}
\eeq
On the other hand, using
\beq
    \det g = \left\{
    \begin{tabular}{l}
      $1 ~~~ (i \ne p)$\\
      $\ep ~~~ (i = p )$ \\
    \end{tabular}
    \right.
    \label{1-25}
\eeq
and (\ref{1-22a}), we find
\beq
    \frac{\d}{\d \ep} F (gZ) =
\left\{
    \begin{tabular}{l}
      $0 ~~~~~~~~~~~~~ (i \ne p)$\\
      $- \frac{1}{\ep^2} F(Z) ~~~ (i = p )$ \\
    \end{tabular}
    \right.
    \label{1-26}
\eeq
Evaluating the derivative at $\ep = 0$ and $\ep = 1$ for
$i \ne p$ and $i = p$, respectively, we then indeed find that
(\ref{1-22a}) leads to the differential equation (\ref{1-2a}).

Similarly, parametrizing $h$ as
\beq
    h = \diag( h_0 , \cdots , h_{j-1} , (1+ \ep )h_j  , h_{j+1} , \cdots , h_n )
    \label{1-27}
\eeq
with $0 \le j \le n$, we find that $Zh$ has only one $\ep$-dependent
column corresponding to the $j$-th column, $\left( z_{0j}( 1 +\ep )h_j  ,
z_{1j} ( 1 +\ep )h_j, \cdots , z_{kj} ( 1 +\ep )h_j \right)^{T}$.
The derivative of $F(Zh)$ with respect to $\ep$ is then expressed as
\beq
    \frac{\d}{\d \ep} F (Zh) \, = \, \sum_{i=0}^{k} z_{ij}
    \frac{\d }{\d z_{ij}} F(Zh) \, = \, \sum_{i=0}^{k} z_{ij}
    \frac{\d }{\d z_{ij}} F(Z)(1+ \ep )^{\al_j} \prod_{ l = 0}^{n} h_{l}^{\al_l}
    \label{1-28}
\eeq
where in the last step we use the relation from (\ref{1-22b}):
\beq
    F (Zh ) \, = \, F(Z) (1+ \ep )^{\al_j} \prod_{ l = 0}^{n} h_{l}^{\al_l}
    \label{1-29}
\eeq
The same derivative can then be expressed as
\beq
    \frac{\d}{\d \ep} F (Zh) \, = \, \al_j F(Z) (1+ \ep )^{\al_j  - 1}
    \prod_{ l = 0}^{k} h_{l}^{\al_l}
    \label{1-30}
\eeq
Setting $\ep = 0$, we can therefore derive the equation (\ref{1-2b}).

The other equation (\ref{1-2c}) for $F(Z)$ follows from the
definition of $\Phi$. From (\ref{1-5a}) and (\ref{1-5b}) we find
that $\Phi$ satisfies
\beq
    \frac{\d \Phi}{\d z_{ip} } \, = \, \frac{ \al_i \tau_p}{l_i (\tau) } \Phi
    \label{1-31}
\eeq
This relation leads to
\beq
    \frac{\d^2 \Phi}{\d z_{ip} \d z_{jq} } \, = \,
    \frac{ \al_i \al_j \tau_p \tau_q}{l_i (\tau) l_j (\tau) } \Phi
    \, = \, \frac{\d^2 \Phi}{\d z_{iq} \d z_{jp} }
    \label{1-32}
\eeq
which automatically derives the equation (\ref{1-2c}).

The integral $F(Z)$ in (\ref{1-4}) therefore indeed satisfies the defining
equations (\ref{1-2a})-(\ref{1-2c}) of the generalized hypergeometric functions on
$Gr ( k+1 , n+1 )$.
{\it The Grassmannian space $Gr ( k+1 , n+1 )$
is defined as a set of $(k+1)$-dimensional linear
subspaces in $(n+1)$-dimensional complex vector space $\C^{n+1}$.}
It is defined as
\beq
    Gr ( k+1 ,n+1 ) \, = \, \widetilde{Z} / GL (k+1 ,\C )
    \label{1-33}
\eeq
where $\widetilde{Z}$ is $(k+1)\times (n+1)$
complex matrices with $rank \widetilde{Z} = k+1$.
Consider some matrix $M$ and assume that there exists
a nonzero $r$-dimensional minor determinant of $M$. Then
the rank of $M$ is in general defined by the largest number of
such $r$'s. Thus $\widetilde{Z}$ is not exactly same as $Z$
defined in (\ref{1-20}). $\widetilde{Z}$ is more relaxed
since it allows some $(k+1)$-dimensional minor determinants vanish,
that is, $Z \subseteq \widetilde{Z}$.
In this sense $F(Z)$ is conventionally called
the generalized hypergeometric functions on $Gr (k+1,n+1)$
and we follow this convention in the present chapter.

\subsection{Non-projected formulation}

In terms of the homogeneous coordinate
$\tau = ( \tau_0 , \tau_1 , \cdots , \tau_k )$ on $\cp^k$,
coordinates on $\C^{k}$ can be parametrized as
\beq
    t_1 = \frac{\tau_1}{\tau_0} ~, ~ t_1 = \frac{\tau_1}{\tau_0} ~,~
    \cdots ~, ~ t_k = \frac{\tau_k}{\tau_0}
    \label{1-34}
\eeq
For simplicity, we now fix $(z_{00} , z_{10} , \cdots , z_{n0} )^T$ at
$(1, 0, \cdots, 0)^T$, {\it i.e.},
\beq
    Z = \left(
          \begin{array}{ccccc}
            1 & z_{01} & z_{02} & \cdots & z_{0n} \\
            0 & z_{11} & z_{12} & \cdots & z_{1n} \\
            \vdots & \vdots & \vdots &  & \vdots \\
            0 & z_{k1} & z_{k2} & \cdots & z_{kn} \\
          \end{array}
        \right)
    \label{1-35}
\eeq
Then the integrand of $F(Z)$ can be expressed as
\beqar
    \Phi \om & = &
    \tau_{0}^{\al_0}  \prod_{j = 1}^{n}
    \left( \tau_0 z_{0j} + \tau_1 z_{1j} + \cdots + \tau_k z_{kj} \right)^{\al_j} \,
    \nonumber \\
    && ~~ \times
    \sum_{i = 0}^{k} ( -1 )^i \tau_i
    d \tau_0 \wedge d \tau_1 \wedge \cdots \wedge
    d \tau_{i-1} \wedge d \tau_{i + 1} \wedge \cdots \wedge d \tau_k
    \nonumber \\
    &=&
    \prod_{j = 1}^{n}
    \left(  z_{0j} + \frac{\tau_1}{\tau_0} z_{1j} + \cdots + \frac{\tau_k}{\tau_0} z_{kj} \right)^{\al_j}
    d \left( \frac{\tau_1}{\tau_0} \right) \wedge d \left( \frac{\tau_2}{\tau_0} \right)
    \wedge \cdots \wedge d \left( \frac{\tau_k}{\tau_0} \right)
    \nonumber \\
    & =& \widetilde{\Phi} \widetilde{\om}
    \label{1-36}
\eeqar
where we use (\ref{1-3b}) and define $\widetilde{\Phi}$, $\widetilde{\om}$ by
\beqar
    \widetilde{\Phi}  &=& \prod_{j=1}^{n} \widetilde{l}_j (t)^{\al_j}
    \label{1-37a} \\
    \widetilde{l}_j (t) &=& z_{0j} + t_1 z_{1j} + t_2 z_{2j} + \cdots + t_k z_{kj}
    ~~~~~ ( 1 \le j \le n )
    \label{1-37b} \\
    \widetilde{\om} &=& dt_1 \wedge dt_2 \wedge \cdots \wedge dt_k
    \label{1-37c}
\eeqar
The exponents $\al_j$ $(j= 1,2,\cdots, n)$ are also imposed to
the non-integer conditions $\al_j \not\in  \Z $
and $\al_1 + \al_2 + \cdots + \al_n \not\in  \Z$.
The multivalued function $\widetilde{\Phi}$ is now defined
in the following space
\beq
    \widetilde{X} \, = \, \C^k - \bigcup_{j= 1}^{n} \widetilde{{\cal H}}_j
    \label{1-38}
\eeq
where
\beq
    \widetilde{{\cal H}}_j \, = \,  \{ t \in \C^k \, ; ~ \widetilde{l}_j (t) = 0 \}
    \label{1-39}
\eeq
These are non-projected versions of (\ref{1-6}) and (\ref{1-7}).

As before, from $\widetilde{\Phi}$ we can define
rank-1 local systems $\widetilde{\cal L}$,
$\widetilde{\cal L}^{\vee}$ on $\widetilde{X}$, which lead
to the $k$-th homology and cohomology groups,
$H_k ( \widetilde{X}, \widetilde{\cal L}^{\vee})$
and $H^k ( \widetilde{X}, \widetilde{\cal L})$.
Then the integral over $\widetilde{\Phi} \widetilde{\om}$ is defined as
\beq
    F(Z) \, = \, \int_{\widetilde{\Del}}\widetilde{\Phi} \widetilde{\om}
    \label{1-40}
\eeq
where $[ \widetilde{\Del} ] = H_k ( \widetilde{X}, \widetilde{\cal L}^{\vee})$
and $[ \widetilde{\om} ] = H^k ( \widetilde{X}, \widetilde{\cal L})$.

In regard to the cohomology group $H^k ( \widetilde{X}, \widetilde{\cal L})$,
Aomoto shows the following theorem\footnote{Theorem 9.6.2 in \cite{Aomoto:1994bk}}:
\begin{enumerate}
  \item The dimension of $H^k ( \widetilde{X}, \widetilde{\cal L})$ is given
  by $\left(
       \begin{array}{c}
        \!\! n-1 \!\! \\
        \!\! k \!\! \\
       \end{array}
     \right)$.
  \item The basis of $H^k ( \widetilde{X}, \widetilde{\cal L})$ can be formed by
  $d \log \widetilde{l}_{j_1} \wedge d \log \widetilde{l}_{j_1} \wedge
  \cdots \wedge d \log \widetilde{l}_{j_k}$ where
  $1\le j_1 < j_2 < \cdots  < j_k \le n-1$.
\end{enumerate}
Correspondingly, the homology group
$H_k ( \widetilde{X}, \widetilde{\cal L}^{\vee})$ has dimension
$\left(
   \begin{array}{c}
    \!\! n-1 \!\! \\
    \!\! k \!\! \\
   \end{array}
\right)$ and its basis can be formed
finite regions bounded by $\widetilde{\cal H}_j$.
In terms of $\widetilde{l}_j$'s the basis
of $H^k ( \widetilde{X}, \widetilde{\cal L})$ can also be
chosen as \cite{Haraoka:2002bk}:
\beq
    \varphi_{j_1 j_2 \dots j_k} \, = \,
    d \log \frac{ \widetilde{l}_{j_1 + 1}}{ \widetilde{l}_{j_1}}
    \wedge
    d \log \frac{ \widetilde{l}_{j_2 + 1}}{ \widetilde{l}_{j_2}}
    \wedge
    \cdots
    \wedge
    d \log \frac{ \widetilde{l}_{j_k + 1}}{ \widetilde{l}_{j_k}}
    \label{1-41}
\eeq
where $1\le j_1 < j_2 < \cdots  < j_k \le n-1$.

\section{Generalized hypergeometric functions on $Gr ( 2, n + 1 )$}

In this section we consider a particular case of $Gr( 2, n+1)$.
The corresponding configuration space is simply given
by $n+1$ distinct points in $\cp^1$. This can be represented
by a $2 \times (n+1)$ matrix $Z$ any of whose 2-dimensional
minor determinants are nonzero.
Allowing the freedom of the
coordinate transformations $GL(2, \C)$ from the left
and the scale transformations $H_2 = \diag(h_0, h_1)$ from the right, we can
uniquely parametrize $Z$ as
\beq
    Z =
    \left(
      \begin{array}{cccccc}
        1 & 0 & 1 & 1 & \cdots & 1 \\
        0 & 1 & -1 & -z_3 & \cdots & -z_{n} \\
      \end{array}
    \right)
    \label{2-1}
\eeq
where $z_i \ne 0 , 1, z_j$ ($i\ne j$, $3 \le i, j \le n$).
Thus we can regard $Z$ as
\beq
    Z \, \simeq \, \{ (z_3 , z_4 , \cdots , z_{n} ) \in \C^{n-2} \, | \,
    z_i \ne 0, 1, z_j ~ (i\ne j) \}
    \label{2-2}
\eeq
The three other points $(z_0 , z_1 , z_2 )$ can be fixed at $\{ 0, 1, \infty \}$
This agrees with the fact that the $GL(2, \C)$
invariance fixes three points out of the $(n+1)$ distinct points in $\cp^1$.

In application of the previous section, we can carry out
a systematic formulation of the generalized hypergeometric functions
on $Gr(2, n+1)$ as follows.
We begin with a multivalued function of a form
\beq
    \Phi =   1^{\al_0} \cdot t^{\al_1} ( 1 -t)^{\al_2}
    (1 - z_3 t)^{\al_3} \cdots ( 1 - z_{n} t )^{\al_n} \, = \, \prod_{j=1}^{n} l_j (t)^{\al_j}
    \label{2-3}
\eeq
where
\beq
    l_0(t) = 1 \, , ~ l_1(t) = t \, , ~ l_2 (t) = 1 -t
    \, , ~ l_j (t) = 1 - z_{j} t ~~ (3 \le j \le n)
    \label{2-4}
\eeq
As in (\ref{1-3a}) and (\ref{1-3b}), the exponents obey the non-integer conditions
\beq
    \al_j \not\in \Z ~~ (0 \le j \le n)\, ,~~~  \sum_{j=0}^{n} \al_j = -2
    \label{2-5}
\eeq
As considered before, the latter condition applies to the expressions
(\ref{1-37c})-(\ref{1-37c}),
that is, when $F(Z)$ is expressed as $F(Z) = \int_\Del \Phi dt$.
The defining space of $\Phi$ is given by
\beq
    X = \cp^1 - \{ 0 ,1, 1/ z_{3} , \cdots , 1/ z_{n} , \infty \}
    \label{2-6}
\eeq
From $\Phi$ we can determine a rank-1 local system ${\cal L}$ on $X$
and its dual local system ${\cal L}^{\vee}$.
Applying the result in (\ref{1-41}), the basis of the cohomology
group $H^1 ( X , {\cal L})$ is then given by
\beq
    d \log \frac{l_{j+1}}{l_j} ~~~~~~ (0 \le j \le n-1)
    \label{2-7}
\eeq
In the present case the basis of the homology group $H_1 ( X , {\cal L}^{\vee} )$
can be specified by a set of paths connecting the branch points.
For example, we can choose these by
\beq
    \Del_{\infty 0} \, , \,
    \Del_{01} \, , \, \Del_{1\frac{1}{z_3}} \, , \, \Del_{\frac{1}{z_3}\frac{1}{z_4} }
    \, , \, \cdots \, , \, \Del_{\frac{1}{z_{n-1}}\frac{1}{z_{n}} }
    \label{2-8}
\eeq
where $\Del_{pq}$ denotes a path on $\cp^1$ connecting branch points $p$ and $q$.
To summarize, for an element $\Del \in H_1 ( X , {\cal L}^{\vee} )$
associated with $\Phi$ of (\ref{2-3}),  we can define
a set of generalized hypergeometric functions on $Gr ( 2, n+1)$ as
\beq
    f_j ( Z ) \, = \, \int_\Del \Phi \, d \log \frac{l_{j+1}}{l_j}
    \label{2-9}
\eeq
where $0 \le j \le n-1$.
In the next section we consider the case of $n=3$, the simplest
case where only one variable exists, which corresponds to Gauss' hypergeometric function.

\section{Reduction to Gauss' hypergeometric function}
\subsection{Basics of Gauss' hypergeometric function}

We first review the basics of Gauss' hypergeometric function.
In power series, it is defined as
\beq
    F (a, b, c; z) \, = \, \sum_{n= 1}^{\infty}
    \frac{ (a)_n (b)_n }{ (c)_n \, n! } z^n
    \label{2-10}
\eeq
where $|z| < 1$, $c \not \in \Z_{\le 0}$ and
\beq
    (a)_{n} =
    \left\{
      \begin{array}{ll}
        1 & ~~ (n=1) \\
        a(a+1)(a+2) \cdots (a+n-1) & ~~(n \ge 1) \\
      \end{array}
    \right.
    \label{2-11}
\eeq
$F (a, b, c; z)$ satisfies the hypergeometric differential equation
\beq
    \left[
    \frac{d^2}{dz^2} + \left( \frac{c}{z} + \frac{a+b+1 - c}{z-1} \right) \frac{d}{d z}
    + \frac{ab}{z (z-1) }
    \right] F ( a, b,c; z)  \, = \, 0
    \label{2-12}
\eeq
Euler's integral formula for $F ( a, b,c; z)$ is written as
\beq
    F ( a, b,c; z) \, = \, \frac{\Ga ( c) }{\Ga( a ) \Ga (c - a)}
    \, \int_{0}^{1} t^{a-1} (1 -t )^{c-a-1} (1 - z t )^{-b} \, dt
    \label{2-13}
\eeq
where $| z | < 1 $ and $0 < \Re(a) < \Re(c ) $\footnote{
This condition can be relaxed to $a \not \in \Z$,
$c-a \not \in \Z$ by use of the well-known Pochhammer contour
in the integral (\ref{2-13}).
}.
$\Ga (a)$'s denote the Gamma functions
\beq
    \Gamma (a ) \, = \, \int_{0}^{\infty} e^{-t} t^{a -1} dt
    ~~~~~ ( \Re (a) > 0 )
    \label{2-13a}
\eeq
The second order differential equation (\ref{2-12}) has regular singularities
at $z = 0, 1, \infty$.
Two independent solutions around each singular point are expressed as
\beqar
    z = 0 ~: && \left\{
    \begin{array}{l}
      f_1 (z) \, = \, F(a,b,c; z) \\
      f_2 (z) \, = \, z^{1-c} F( a-c+1, b-c+1, 2-c;z)
    \end{array}
    \right.
    \label{2-14a} \\
    z = 1  ~: && \left\{
    \begin{array}{l}
      f_3 (z) \, = \, F(a,b,a+b-c+1; 1-z) \\
      f_4 (z) \, = \, (1-z)^{c-a-b} F( c-a, c-a, c-a-b+1;1-z)
    \end{array}
    \right.
    \label{2-14b} \\
    z = \infty ~: &&  \left\{
    \begin{array}{l}
      f_5 (z) \, = \, z^{-a} F(a,a-c+1,a-b+1; 1/z) \\
      f_6 (z) \, = \, z^{-b} F(b-c+1, b,b-a+1;1z)
    \end{array}
    \right.
    \label{2-14c}
\eeqar
where we assume $c \not \in \Z$,
$a+b-c \not \in \Z$
and $a-b \not \in \Z$
at $z= 0$, $z = 1$
and $z = \infty$, respectively.

\subsection{Reduction to Gauss' hypergeometric function 1: From defining equations}

From (\ref{2-13}) we find the relevant $2 \times 4$ matrix in a form of
\beq
    Z =
    \left(
      \begin{array}{cccc}
        1 & 0 & 1 & 1  \\
        0 & 1 & -1 & -z  \\
      \end{array}
    \right)
    \label{2-15}
\eeq
The set of equations (\ref{1-2a})-(\ref{1-2c})
then reduce to the followings:
\beqar
    ( \d_{00} + \d_{02} + \d_{03} ) F &=& -F
    \label{2-16a} \\
    ( \d_{11} - \d_{12} + z \d_z ) F &=& -F
    \label{2-16b} \\
    ( \d_{10} + \d_{12} - \d_{z} ) F &=& 0
    \label{2-16c} \\
    ( \d_{01} + \d_{02} - \d_{03} ) F &=& 0
    \label{2-16d} \\
    \d_{00} F &=& \al_0 F
    \label{2-17a} \\
    \d_{11} F &=& \al_1 F
    \label{2-17b} \\
    ( \d_{02} - \d_{12} ) F &=& \al_2 F
    \label{2-17c} \\
    ( \d_{03} + z \d_{z} ) F &=& \al_3 F
    \label{2-17d} \\
    -  \d_{z} \d_{02} F &= & \d_{12} \d_{03} F
    \label{2-18}
\eeqar
where $\d_{ij} = \frac{\d}{\d z_{ij}}$ and $\d_{13} = - \frac{\d}{\d z} = - \d_z$.
The last relation (\ref{2-18}) arises from (\ref{1-2c}); we here write
down the one that is nontrivial and involves $\d_z$.
Since the sum of (\ref{2-16a}) and (\ref{2-16b}) equals to
the sum of (\ref{2-17a})-(\ref{2-17d}), we can easily find
$\al_0 + \cdots + \al_3 = -2$ in accord with (\ref{2-5}).
The second order equation (\ref{2-18}) is then expressed as
\beq
    - \d_z ( \al_1 + \al_2 + 1 + z \d_z  ) F \, = \,
    ( \al_1 + 1 + z \d_z ) ( \al_3 - z \d_z ) F
    \label{2-19}
\eeq
This can also be written as
\beq
    \left[
    z (1-z) \d^{2}_{z} + \left( c - (a+ b+1) z \right) \d_z -ab
    \right] F \, = \, 0
    \label{2-20}
\eeq
where
\beqar
    a &=& \al_1 + 1
    \nonumber \\
    b &=& - \al_3
    \label{2-21} \\
    c &=& \al_1 + \al_2 + 2
    \nonumber
\eeqar
We can easily check that (\ref{2-20})
identifies with the hypergeometric differential equation (\ref{2-12}).

As seen in (\ref{2-1}), there exist multiple complex variables
for $n > 3$. In these cases reduction of the defining equations
(\ref{1-2a})-(\ref{1-2c}) can be carried out in principle but, unfortunately,
is not as straightforward as the case of $n=3$.

\subsection{Reduction to Gauss' hypergeometric function 2: Use of twisted cohomology}

The hypergeometric equation (\ref{2-12}) is a second order differential equation.
Setting $f_1 = F$, $f_2 =  \frac{z}{b} \frac{d}{dz} F$, we can express (\ref{2-12})
in a form of a first order Fuchsian differential equation \cite{Aomoto:1994bk}:
\beq
    \frac{d}{d z}
    \left(
      \begin{array}{c}
        f_1 \\
        f_2 \\
      \end{array}
    \right)
    =
    \left(
    \frac{A_0}{z} + \frac{A_1}{z-1}
    \right)
     \left(
      \begin{array}{c}
        f_1 \\
        f_2 \\
      \end{array}
    \right)
    \label{2-22}
\eeq
where
\beq
    A_0 =
    \left(
      \begin{array}{cc}
        0 & b \\
        0 & 1-c \\
      \end{array}
    \right),
    ~~
    A_1 =
    \left(
      \begin{array}{cc}
        0 & 0 \\
        -a & c-a-b-1 \\
      \end{array}
    \right)
    \label{2-23}
\eeq
Using the results (\ref{2-3})-(\ref{2-9}), we now
obtain other first order representations of
the hypergeometric differential equation.

Let us start with a {\it non-projected} multivalued function
\beq
    \Phi \, = \, t^a (1-t)^{c-a} ( 1- zt )^{-b}
    \label{2-24}
\eeq
where
\beq
    a \, , \, c-a \, , \, -b \not \in \Z
    \label{2-25}
\eeq
$\Phi$ is defined on $X = \cp^1 - \{ 0,1,1/z , \infty \}$.
From these we can determine
a rank-1 local system ${\cal L}$ and its dual ${\cal L}^\vee$ on $X$.
Then, using (\ref{2-7}), we can obtain
a basis of the cohomology group $H^1 ( X , {\cal L})$ given by
the following set
\beqar
    \varphi_{\infty 0} &=& \frac{dt}{t}
    \label{2-26a} \\
    \varphi_{01} &=& \frac{dt}{t (1-t)}
    \label{2-26b} \\
    \varphi_{1 \frac{1}{z}} &=& \frac{(z-1) dt}{(1-t) ( 1 -zt)}
    \label{2-26c}
\eeqar
Similarly, from (\ref{2-8}) a basis of the homology group $H_1 ( X ,{\cal L}^\vee )$
is given by
\beq
    \{ \Del_{\infty 0} \, , \, \Del_{01} \, , \, \Del_{1 \frac{1}{z}}
    \}
    \label{2-27}
\eeq
In terms of these we can express Gauss' hypergeometric function as
\beq
    f_{01} (Z)
    \,  =  \,  \int_{\Del_{01}}  \Phi \varphi_{01}
    \, = \,
    \int_{0}^{1} t^{a-1} ( 1- t)^{c-a-1} ( 1 - z t)^{-b}  dt
    \label{2-28}
\eeq

The derivative of $f_{01} (Z) = f_{01}(z)$ with respect to $z$ is written as
\beq
    \frac{d}{d z} f_{01}(z) \, = \,
    \frac{d}{d z} \int_{\Del_{01}}  \Phi \varphi_{01}
    \, = \,
    \int_{\Del_{01}}  \Phi  \, \nabla_{\! z} \varphi_{01}
    \label{2-29}
\eeq
where
\beq
    \nabla_{\! z} \, = \, \d_z + \d_z \log \Phi
    \, = \, \d_z + \frac{bt}{1-zt}
    \label{2-30}
\eeq
Thus the derivative comes down to the computation of
$\nabla_{\! z} \varphi_{01}$; notice that the choice of a twisted cycle $\Del$
is irrelevant as far as the derivative itself is concerned.
In order to make sense of (\ref{2-29}) we should require
$\nabla_{\! z} \varphi_{01} \in H^1 ( X , {\cal L})$, that
is, it should be represented by a linear combinations of
(\ref{2-26a})-(\ref{2-26c}).
There is a caveat here, however. We know that
an element of $H^1 ( X , {\cal L})$ forms an equivalent class
as discussed earlier; see (\ref{1-8}) and (\ref{1-9}).
In the present case ($k=1$), $\al$ in (\ref{1-8}) is a 0-form
or a constant. So we can demand
\beq
    d \log \Phi
    \, = \,
    a \frac{dt}{t} - (c-a) \frac{dt}{1-t}
    + b \frac{z \, dt}{1-zt}
    \, \equiv \,
    0
    \label{2-31}
\eeq
in the computation of $\nabla_{\! z} \varphi_{01}$.
This means that the number of the base elements
can be reduced from 3 to 2. Namely, any elements of
$H^1 ( X , {\cal L})$ can be expressed by a combinations
of an arbitrary pair in (\ref{2-26a})-(\ref{2-26c})
under the condition (\ref{2-31}).
This explains the
numbering discrepancies between (\ref{1-41}) and
(\ref{2-7}) and agrees with the general result in the previous section that
the dimension of the cohomology group is given by
$\left(
   \begin{array}{c}
    \!\! n-1 \!\! \\
     \!\! k \!\! \\
   \end{array}
 \right) =
 \left(
   \begin{array}{c}
    \!\! 2 \!\! \\
    \!\! 1 \!\! \\
   \end{array}
 \right) = 2$.

Choosing the pair of $( \varphi_{01}, \varphi_{\infty 0} )$, we find
\beqar
    \nabla_{\! z} \varphi_{\infty 0} &=& \frac{b dt}{1-zt}
    \nonumber \\
    & \equiv &
    \frac{1}{z} \left( - a \frac{dt}{t} + (c-a) \frac{dt}{1-t} \right)
    \nonumber \\
    &=& \frac{c-a}{z} \varphi_{01} - \frac{c}{z} \varphi_{\infty 0}
    \label{2-32a} \\
    \nabla_{\! z} \varphi_{01} &=&
    \nabla_{\! z} \left( \varphi_{\infty 0} + \frac{dt}{1-t} \right)
    \nonumber \\
    &=&
    \nabla_{\! z} \varphi_{\infty 0} + \frac{b}{1-z} \left( \frac{dt}{1-t} - \frac{dt}{1-zt} \right)
    \nonumber \\
    &=&
    \frac{z}{z-1} \nabla_{\! z} \varphi_{\infty 0} - \frac{b}{z-1} ( \varphi_{01} - \varphi_{\infty 0} )
    \nonumber \\
    & \equiv &
    \frac{c-a-b}{z-1} \varphi_{01} + \frac{b-c}{z-1} \varphi_{\infty 0}
    \label{2-32b}
\eeqar
where notation $\equiv$ means the use of condition (\ref{2-31}).
Using (\ref{2-29}), we obtain a first order differential equation
\beq
    \frac{d}{dz}
    \left(
       \begin{array}{c}
         f_{01} \\
         f_{\infty 0} \\
       \end{array}
    \right)
    =
    \left(
        \frac{A^{(\infty 0)}_{0}}{z}   + \frac{A^{(\infty 0)}_{1}}{z -1 }
    \right)
    \left(
       \begin{array}{c}
         f_{01} \\
         f_{\infty 0} \\
       \end{array}
     \right)
    \label{2-33}
\eeq
where
\beq
    A^{(\infty 0)}_{0} =
    \left(
      \begin{array}{cc}
        0 & 0 \\
        c-a & -c \\
      \end{array}
    \right),
    ~~
    A^{(\infty 0)}_{1} =
    \left(
      \begin{array}{cc}
        c-a-b & b-c \\
        0 & 0 \\
      \end{array}
    \right)
    \label{2-34}
\eeq
Solving for $f_{01}$, we can easily confirm that (\ref{2-33}) leads
to Gauss' hypergeometric differential equation (\ref{2-12}).

Similarly, for the choice of
$( \varphi_{01}, \varphi_{1 \frac{1}{z} } )$ we find
\beqar
    \nabla_{\! z} \varphi_{01}
    &=&
    \frac{b}{z-1} \varphi_{1 \frac{1}{z} }
    \label{2-35a} \\
    \nabla_{\! z} \varphi_{1 \frac{1}{z} }
    & \equiv &
    \nabla_{\! z} \left( -\frac{a}{b} \frac{z-1}{z} \varphi_{01} +
    \frac{c-a}{b} \frac{z-1}{z} \frac{dt }{(1-t)^2} \right)
    \nonumber \\
    & \equiv &
    - \frac{a}{z} \varphi_{01} + \left( - \frac{c+1}{z} + \frac{c-a-b+1}{z-1} \right)
    \varphi_{1 \frac{1}{z} }
    \label{2-35b}
\eeqar
Notice that
$\varphi_{01}$ and $\varphi_{1 \frac{1}{z} }$ have the same factor
$(1-t)^{-1}$. This factor can be absorbed in the definition of $\Phi$ in (\ref{2-24}).
Thus, in applying the derivative formula (\ref{2-29}), we should replace $c$
by $c-1$. This leads to another first order differential equation
\beq
    \frac{d}{dz}
    \left(
       \begin{array}{c}
         f_{01} \\
         f_{1 \frac{1}{z}} \\
       \end{array}
    \right)
    =
    \left(
        \frac{A^{(1 \frac{1}{z})}_{0}}{z}   + \frac{A^{(1 \frac{1}{z})}_{1}}{z -1 }
    \right)
    \left(
       \begin{array}{c}
         f_{01} \\
         f_{1 \frac{1}{z}} \\
       \end{array}
     \right)
    \label{2-36}
\eeq
where
\beq
    A^{(1 \frac{1}{z})}_{0} =
    \left(
      \begin{array}{cc}
        0 & 0 \\
        -a & -c \\
      \end{array}
    \right),
    ~~
    A^{(1 \frac{1}{z})}_{1} =
    \left(
      \begin{array}{cc}
        0 & b \\
        0 & c-a-b \\
      \end{array}
    \right)
    \label{2-37}
\eeq
Solving for $f_{01}$, we can also check that
(\ref{2-36}) becomes Gauss' hypergeometric
differential equation (\ref{2-12}).

The representations (\ref{2-23}) and (\ref{2-34})
are obtained by Aomoto-Kita \cite{Aomoto:1994bk} and Haraoka
\cite{Haraoka:2002bk}, respectively.
The last one (\ref{2-37}) is not known in the literature
as far as the author notices.
Along the lines of the above derivation, we can also obtain
the Aomoto-Kita representation (\ref{2-23}) as follows.
We introduce a new one-form
\beq
    \widetilde{\varphi}_{1 \frac{1}{z}}
    \, =  \, \frac{z}{z-1} \varphi_{1 \frac{1}{z}}
    \, = \,  \frac{z \, dt}{ (1-t)(1-zt) }
    \label{2-38}
\eeq
The corresponding hypergeometric function is given by
$\widetilde{f}_{1 \frac{1}{z}} = \int_{\Del_{01}}
\Phi \widetilde{\varphi}_{1 \frac{1}{z}}$.
From (\ref{2-35a}) we can easily see $\nabla_{\! z} \varphi_{01} =
\frac{b}{z} \widetilde{\varphi}_{1 \frac{1}{z}}$. This is
consistent with the condition
$f_1 = F$, $f_2 =  \frac{z}{b} \frac{d}{dz} F$ in (\ref{2-22}).
Since $z$ is defined as $z \ne 0 ,1$,
$\frac{b}{z} \widetilde{\varphi}_{1 \frac{1}{z}}$ and
$\frac{b}{z-1} \varphi_{1 \frac{1}{z}}$ are equally well defined one-forms.
We can then choose
the pair $( \varphi_{01} , \widetilde{\varphi}_{1 \frac{1}{z}} )$
as a possible basis of the cohomology group.
The derivatives $\nabla_{\! z} \varphi_{01}$,
$\nabla_{\! z} \widetilde{\varphi}_{1 \frac{1}{z}} $ are calculated as
\beqar
    \nabla_{\! z} \varphi_{01}
    &=&
    \frac{b}{z} \widetilde{\varphi}_{1 \frac{1}{z}}
    \label{2-39a} \\
    \nabla_{\! z} \widetilde{\varphi}_{1 \frac{1}{z}}
    & \equiv &
    \nabla_{\! z} \left(
    - \frac{a}{b} \varphi_{01} + \frac{c-a}{b} \frac{dt}{ (1-t)^2 }
    \right)
    \nonumber \\
    & \equiv &
    - \frac{a}{z-1}  \varphi_{01} + \left( \frac{-c}{z} + \frac{c-a-b}{z-1}
    \right)\widetilde{\varphi}_{1 \frac{1}{z}}
    \label{2-39b}
\eeqar
where we use the relations
\beqar
    \frac{t \, dt }{(1-zt)( 1-t)} &=& \frac{1}{z-1} \left(
    \frac{dt}{1-zt} - \frac{dt}{1-t} \right)
    \label{2-40a} \\
    \frac{dt}{(1-t)^2} & \equiv & \frac{1}{c-a}
    \left( a \varphi_{01} + b \widetilde{\varphi}_{1 \frac{1}{z}} \right)
    \label{2-40b}
\eeqar
As before, $\varphi_{01}$ and $\widetilde{\varphi}_{1 \frac{1}{z}}$ have
the same factor $(1-t)^{-1}$. Thus, replacing $c$ by $c-1$, we obtain
a first order differential equation
\beq
    \frac{d}{dz}
    \left(
       \begin{array}{c}
         f_{01} \\
         \widetilde{f}_{1 \frac{1}{z}} \\
       \end{array}
    \right)
    =
    \left(
        \frac{\widetilde{A}^{(1 \frac{1}{z})}_{0}}{z}   + \frac{\widetilde{A}^{(1 \frac{1}{z})}_{1}}{z -1 }
    \right)
    \left(
       \begin{array}{c}
         f_{01} \\
         \widetilde{f}_{1 \frac{1}{z}} \\
       \end{array}
     \right)
    \label{2-41}
\eeq
where
$\widetilde{f}_{1 \frac{1}{z}} = \int_{\Del_{01}} \Phi \widetilde{\varphi}_{1 \frac{1}{z}}$
and
\beq
    \widetilde{A}^{(1 \frac{1}{z})}_{0} =
    \left(
      \begin{array}{cc}
        0 & b \\
        0 & 1-c \\
      \end{array}
    \right),
    ~~
    \widetilde{A}^{(1 \frac{1}{z})}_{1} =
    \left(
      \begin{array}{cc}
        0 & 0 \\
        -a & c-a-b-1 \\
      \end{array}
    \right)
    \label{2-42}
\eeq
We therefore reproduce the Aomoto-Kita representation (\ref{2-22}), (\ref{2-23})
by a systematic construction of first order representations of the hypergeometric
differential equation.

Lastly, we note that $\varphi_{\infty 0} = \frac{dt}{t}$ and $\varphi_{01} = \frac{dt}{t(1-t)}$
have the same factor $t^{-1}$ but we can not absorb this factor into $\Phi$.
This is because we can not obtain $dt$ as a base element of
$H^1 ( X , {\cal L})$ which is generically given in a form of
$d \log \frac{l_{j+1}}{l_j}$ as discussed in (\ref{2-7}).

\subsection{Reduction to Gauss' hypergeometric function 3: Permutation invariance}

The choice of twisted cycles or $\Del$'s is irrelevant in the above
derivations of the first order Fuchsian differential equations.
The hypergeometric function is therefore
satisfied by a more general integral form, rather than (\ref{2-28}), {\it i.e.},
\beq
    f_{01}^{(\Del_{pq})} (z) \, = \,
    \int_{\Del_{pq}} \Phi \varphi_{01}
    \, = \, \int_{p}^{q} t^{a-1} ( 1- t)^{c-a-1} ( 1 - z t)^{-b}
    dt
    \label{2-43}
\eeq
where $(p, q)$ represents an arbitrary pair among the four branch points
$p, q \in \{ 0, 1, 1/z , \infty \}$.
This means that we can impose {\it permutation invariance} on
the branch points. $\Del_{pq}$ is then given by the following set of
twisted cycles:
\beq
    \Del_{pq} \, = \,
    \{  \Del_{\infty 0} \, , \, \Del_{01} \, , \, \Del_{1 \frac{1}{z}}
    \, , \,  \Del_{1 \infty} \, , \, \Del_{\frac{1}{z} \infty} \, , \, \Del_{0 \frac{1}{z}}
    \}
    \label{2-44}
\eeq
so that the number of elements becomes
$\left(
   \begin{array}{c}
    \!\! 4 \!\! \\
    \!\! 2 \!\! \\
   \end{array}
\right) = 6$. Correspondingly, the base elements of the cohomology group also include
\beqar
    \varphi_{1 \infty } &=& \frac{dt}{1-t}
    \label{2-45a} \\
    \varphi_{\frac{1}{z} \infty} &=& \frac{z \, dt}{1-zt}
    \label{2-45b} \\
    \varphi_{0 \frac{1}{z}} &=& \frac{ dt}{t ( 1 -zt)}
    \label{2-45c}
\eeqar
besides (\ref{2-26a})-(\ref{2-26c}).
It is known that $f_{01}^{(\Del_{pq})} (z)$ are related to
the local solutions $f_i (z)$ $(i= 1,2, \cdots , 6)$ in (\ref{2-14a})-(\ref{2-14c}) by
\beqar
      f_{01}^{(\Del_{01})} (z) & = & B( a, c-a) f_1 (z)
      \label{2-46a} \\
      f_{01}^{ ( \Del_{\frac{1}{z} \infty} )} (z) & = & e^{i \pi (a+b-c+1)} B( b-c+1 , 1-b) f_2 (z)
      \label{2-46b} \\
      f_{01}^{(\Del_{\infty 0})} (z) & = & e^{i \pi (1-a)} B( a , b-c+1) f_3 (z)
      \label{2-46c} \\
      f_{01}^{(\Del_{1 \frac{1}{z}})} (z) & = & e^{i \pi (a-c+1)} B( c-a , 1-b) f_4 (z)
      \label{2-46d} \\
      f_{01}^{(\Del_{1 \frac{1}{z}})} (z) & = & B( a , 1-b) f_5 (z)
      \label{2-46e} \\
      f_{01}^{(\Del_{1 \infty})} (z) & = & e^{- i \pi (a+b-c+1)} B( b-c+1 , c-a) f_6 (z)
      \label{2-46f}
\eeqar
where $B(a, b)$ is the beta function
\beq
    B (a, b) \, = \, \frac{\Ga (a) \Ga (b)}{\Ga(a+b)}
    \, = \, \int_{0}^{1} t^{a-1} (1-t)^{b-1} dt
    ~~~~~ ( \Re (a ) > 0 \, , \, \Re(b) > 0 )
    \label{2-47}
\eeq
(For derivations and details of these relations, see \cite{Haraoka:2002bk}.)

The relevant configuration space represented by $Z$ is given by
$Gr (2, 4) / \S_4$ where $\S_4$ denotes the rank-4 symmetry group.
The permutation invariance can also be confirmed by
deriving another set of the first order differential equations
with the choice of $\varphi_{01}$ and one of (\ref{2-45a})-(\ref{2-45c}).
This is what we will present in the following.

For the choice of
$( \varphi_{01}, \varphi_{1 \infty} )$ we find
\beqar
    \nabla_{\! z} \varphi_{1 \infty}
    &=&
    \frac{bt}{1-zt} \frac{dt}{(1-t)}
    \nonumber \\
    & \equiv &
    - \frac{a}{z(z-1)} \varphi_{01} + \left( \frac{c}{z(z-1)} - \frac{b}{z-1} \right)
    \varphi_{1 \infty}
    \label{2-48a} \\
    \nabla_{\! z} \varphi_{01 }
    &=& \nabla_{\! z} \frac{dt}{t} + \nabla_{\! z} \varphi_{1 \infty}
    \nonumber \\
    & \equiv &
    - \frac{a}{z-1} \varphi_{01} +  \frac{c-b}{z-1}  \varphi_{1 \infty }
    \label{2-48b}
\eeqar
The corresponding differential equation is then expressed as
\beq
    \frac{d}{dz}
    \left(
       \begin{array}{c}
         f_{01} \\
         f_{1 \infty} \\
       \end{array}
    \right)
    =
    \left(
        \frac{A^{(1 \infty)}_{0}}{z}   + \frac{A^{(1 \infty)}_{1}}{z -1 }
    \right)
    \left(
       \begin{array}{c}
         f_{01} \\
         f_{1 \infty} \\
       \end{array}
     \right)
    \label{2-49}
\eeq
where
\beq
    A^{(1 \infty)}_{0} =
    \left(
      \begin{array}{cc}
        0 & 0 \\
        a & -c \\
      \end{array}
    \right),
    ~~
    A^{(1 \infty)}_{1} =
    \left(
      \begin{array}{cc}
        -a & c-b \\
        -a & c-b \\
      \end{array}
    \right)
    \label{2-50}
\eeq
Solving for $f_{01}$, we can check that
(\ref{2-49}) indeed becomes Gauss' hypergeometric differential equation (\ref{2-12}).

Similarly, for $( \varphi_{01}, \varphi_{\frac{1}{z} \infty} )$ we find
\beqar
    \nabla_{\! z} \varphi_{01 }
    &=& \frac{b \, dt}{(1-zt)(1-t)}
    \nonumber \\
    & \equiv &
    - \frac{1}{z-1} \frac{ab}{c} \varphi_{01} -  \frac{1}{z-1} \frac{b}{c} (b-c)
    \varphi_{\frac{1}{z} \infty }
    \label{2-51a} \\
    \nabla_{\! z} \varphi_{\frac{1}{z} \infty}
    & \equiv &
    \nabla_{\! z} \left( -\frac{a}{b} \varphi_{01} + \frac{c}{b} \frac{dt}{1-t} \right)
    \nonumber \\
    & \equiv &
    \frac{1}{z-1} \frac{a}{c} (a -c) \varphi_{01} + \left( \frac{1}{z-1} \frac{1}{c} (b-c)(a-c)
    - \frac{c}{z} \right)
    \varphi_{\frac{1}{z} \infty}
    \label{2-51b}
\eeqar
The first order differential equation is then expressed as
\beq
    \frac{d}{dz}
    \left(
       \begin{array}{c}
         f_{01} \\
         f_{\frac{1}{z} \infty} \\
       \end{array}
    \right)
    =
    \left(
        \frac{A^{(\frac{1}{z} \infty)}_{0}}{z}   + \frac{A^{(\frac{1}{z} \infty)}_{1}}{z -1 }
    \right)
    \left(
       \begin{array}{c}
         f_{01} \\
         f_{\frac{1}{z} \infty} \\
       \end{array}
     \right)
    \label{2-52}
\eeq
where
\beq
    A^{(\frac{1}{z} \infty)}_{0} =
    \left(
      \begin{array}{cc}
        0 & 0 \\
        0 & -c \\
      \end{array}
    \right),
    ~~
    A^{(\frac{1}{z} \infty)}_{1} =
    \left(
      \begin{array}{cc}
        -\frac{ab}{c} & - \frac{b}{c}(b-c) \\
        \frac{a}{c}(a-c) & \frac{1}{c} (b-c)(a-c) \\
      \end{array}
    \right)
    \label{2-53}
\eeq
We can check that (\ref{2-52}) reduces to the hypergeometric
differential equation for $f_{01}$.

Lastly, for $(\varphi_{01} , \varphi_{0 \frac{1}{z}} )$ we find
\beqar
    \nabla_{\! z} \varphi_{01 }
    &=& \frac{b \, dt}{(1-zt)(1-t)}
    \nonumber \\
    & = &
    - \frac{b}{z-1} \left( \varphi_{01} -   \varphi_{0 \frac{1}{z} } \right)
    \label{2-54a} \\
    \nabla_{\! z} \varphi_{0 \frac{1}{z}}
    & \equiv &
    \frac{b \, dt}{ ( 1-zt)^2 }
    \nonumber \\
    & \equiv &
    - \frac{c-a}{z(z-1)}  \varphi_{01} + \frac{c - az}{z(z-1)} \varphi_{0 \frac{1}{z}}
    \label{2-54b}
\eeqar
The corresponding differential equation becomes
\beq
    \frac{d}{dz}
    \left(
       \begin{array}{c}
         f_{01} \\
         f_{0 \frac{1}{z}} \\
       \end{array}
    \right)
    =
    \left(
        \frac{A^{(0 \frac{1}{z})}_{0}}{z}   + \frac{A^{(0 \frac{1}{z})}_{1}}{z -1 }
    \right)
    \left(
       \begin{array}{c}
         f_{01} \\
         f_{0 \frac{1}{z}} \\
       \end{array}
     \right)
    \label{2-55}
\eeq
where
\beq
    A^{(0 \frac{1}{z})}_{0} =
    \left(
      \begin{array}{cc}
        0 & 0 \\
        c-a & -c \\
      \end{array}
    \right),
    ~~
    A^{(0 \frac{1}{z})}_{1} =
    \left(
      \begin{array}{cc}
        -b & b \\
        -c+a & c-a \\
      \end{array}
    \right)
    \label{2-56}
\eeq
We can check that (\ref{2-55}) reduces to the hypergeometric
differential equation for $f_{01}$ as well.

As in the case of (\ref{2-38}), it is tempting to
think of $\widetilde{\varphi}_{ \frac{1}{z} \infty}
=  \frac{z-1}{z} \varphi_{\frac{1}{z} \infty}
=  \frac{(z-1) \, dt}{ 1-zt }$.
But, with $\varphi_{01}$ and $\widetilde{\varphi}_{ \frac{1}{z} \infty}$,
it is not feasible to obtain a first order
differential equation in the form of (\ref{2-52})
which leads to the hypergeometric differential equation.
This is because, if expanded in $\varphi_{01}$ and $\widetilde{\varphi}_{ \frac{1}{z} \infty}$,
the $z$-dependence of the derivatives
$\nabla_{\! z}  \varphi_{01}$ and $\nabla_{\! z}  \widetilde{\varphi}_{ \frac{1}{z} \infty}$,
can not be written in terms of $\frac{1}{z}$ or $\frac{1}{z-1}$.

\subsection{Summary}

In this section we carry out a systematic derivation
of first order representations of the hypergeometric differential
equation by use of twisted cohomology as the simplest
reduction of Aomoto's generalized hypergeometric function.
The first order equations are generically expressed as
\beq
    \frac{d}{d z}
    \left(
      \begin{array}{c}
        f_{01} \\
        f_{pq} \\
      \end{array}
    \right)
    =
    \left(
    \frac{A^{(pq)}_{0}}{z} + \frac{A^{(pq)}_{1}}{z-1}
    \right)
     \left(
      \begin{array}{c}
        f_{01} \\
        f_{pq} \\
      \end{array}
    \right)
    = \, A^{(pq) }_{01}
     \left(
      \begin{array}{c}
        f_{01} \\
        f_{pq} \\
      \end{array}
    \right)
    \label{2-57}
\eeq
where $(pq)$ denotes a pair of four branch points
$\{ 0, 1 , 1/z , \infty \}$ in $\Phi = t^a (1-t)^{c-a} (1-zt)^{-b}$.
A list of the $(2 \times 2)$ matrices $A^{(pq)}_{01}$
obtained in this section is given by the following:
\beqar
    A_{01}^{( \infty 0)} &=&
    \left(
      \begin{array}{cc}
        \frac{c-a-b}{z-1} & \frac{b-c}{z-1} \\
        \frac{c-a}{z} & - \frac{c}{z} \\
      \end{array}
    \right)
    \label{2-58a} \\
    A_{01}^{(1 \frac{1}{z})} &=&
    \left(
      \begin{array}{cc}
        0 & \frac{b}{z-1} \\
        - \frac{a}{z} & - \frac{c}{z} + \frac{c-a-b}{z-1} \\
      \end{array}
    \right)
    \label{2-58b} \\
    \widetilde{A}_{01}^{(1 \frac{1}{z})} &=&
    \left(
      \begin{array}{cc}
        0 & \frac{b}{z} \\
        - \frac{a}{z-1} & - \frac{c-1}{z} + \frac{c-a-b-1}{z-1} \\
      \end{array}
    \right)
    \label{2-58c} \\
    A_{01}^{(1 \infty )} &=&
    \left(
      \begin{array}{cc}
        -\frac{a}{z-1} & \frac{c-b}{z-1} \\
        \frac{a}{z} - \frac{a}{z-1} & - \frac{c}{z} + \frac{c-b}{z-1} \\
      \end{array}
    \right)
    \label{2-58d} \\
    A_{01}^{( \frac{1}{z} \infty )} &=&
    \left(
      \begin{array}{cc}
        - \frac{1}{z-1} \frac{ab}{c} & - \frac{1}{z-1} \frac{b}{c} (b-c)  \\
        \frac{1}{z-1} \frac{a}{c} (a-c) & - \frac{c}{z} + \frac{1}{z-1} \frac{1}{c}(b-c) (a-c) \\
      \end{array}
    \right)
    \label{2-58e} \\
    A_{01}^{( 0 \frac{1}{z})} & =&
    \left(
      \begin{array}{cc}
        - \frac{b}{z-1} & \frac{b}{z-1} \\
        \frac{c-a}{z} - \frac{c-a}{z-1} & - \frac{c-1}{z} + \frac{c-a}{z-1} \\
      \end{array}
    \right)
    \label{2-58f}
\eeqar
where we include the Aomoto-Kita representation
$\widetilde{A}_{01}^{(1 \frac{1}{z})}$. As far as the author notices,
these expressions except (\ref{2-58a}, \ref{2-58c}) are new
for the description of the hypergeometric differential equation.
A common feature among these matrices is that the determinant is
identical:
\beq
    \det A^{(pq)}_{01} \, = \, \frac{ab}{z(z-1)}
    \label{2-59}
\eeq
In terms of the first order differential equation (\ref{2-57}),
this means that the action of the derivative on the basis
$
     \left(
      \begin{array}{c}
        f_{01} \\
        f_{pq} \\
      \end{array}
    \right)
$
of the cohomology group $H^1 ( X , {\cal L} )$
can be represented by a generator of the $SL(2, \C)$ algebra.
In other words, a change of the bases
is governed by the $SL(2, \C)$ symmetry.
The $SL(2, \C)$ invariance corresponds to
the global conformal symmetry for holomorphic functions
on $\cp^1$. In the present case we start from
the holomorphic multivalued function $\Phi$ in (\ref{2-24}) which
is defined on $X = \cp^1 - \{ 0 , 1, 1/z , \infty \}$.
The result (\ref{2-59}) is thus natural in concept but nontrivial in practice because
the equivalence condition $d \log \Phi \equiv 0$ in (\ref{2-31})
is implicitly embedded into the expressions (\ref{2-58a})-(\ref{2-58f}).

\vskip 0.5cm

To conclude this chapter, we first review the definition of Aomoto's generalized
hypergeometric functions on $Gr(k+1, n+1)$, interpreting
their integral representations
in terms of twisted homology and cohomology.
We then consider reduction of the general $Gr(k+1 , n+1)$ case to
particular $Gr(2, n+1)$ cases. The case of $Gr( 2, 4)$ leads to
Gauss' hypergeometric functions. We carry out a thorough study
of this case in the present section.
Much of the present chapter, by nature, deals with reviews of
existed literature. But the results in (\ref{2-57})-(\ref{2-59})
are new as far as the author notices.


\end{document}